\documentclass[12pt,reqno]{amsart}

\usepackage{amstext,amsmath,amssymb,amsthm}
\usepackage[dvips]{graphics,color}

\numberwithin{equation}{section}

\renewcommand{\rq}[1]{(\ref{#1})}
\newcommand{\dsize}{\displaystyle}

\renewcommand{\Re}{\mbox{\rm Re}}

\newcommand{\R}{{ \Bbb  R  }}
\newcommand{\C}{ \Bbb  C }
\newcommand{\N}{ \Bbb  N }

\newcommand{\BEL}{\begin{equation}\label}
\newcommand{\EE}{\end{equation}}

\renewcommand{\medskip}{\vskip .5 cm}
\newtheorem{Thm}{Theorem}[section]
\newtheorem{Lemma}[Thm]{Lemma}
\newtheorem{Cor}[Thm]{Corollary}

\newtheorem{dfn}[Thm]{Definition}
\newtheorem{Ex}[Thm]{Example}

\begin{document}

\title[Minimal potential results in a slab]{Minimal potential results for the Schr\"odinger equation in a slab}

\author[De Carli]{Laura De Carli}
\address{Department of Mathematics and Statistics, Florida International University, Miami,
FL 33199, U.S.A.} \email{decarlil@fiu.edu}

\author[Hudson]{Steve Hudson}
\address{Department of Mathematics and Statistics, Florida International University, Miami,
FL 33199, U.S.A.} \email{hudsons@fiu.edu}

\author[Li]{Xiaosheng Li}
\address{Department of Mathematics and Statistics, Florida International University, Miami,
FL 33199, U.S.A.} \email{xli@fiu.edu}
\thanks{X. Li is supported in part by NSF grant DMS-1109561.}

 \maketitle

\begin{abstract}
Consider the Schr\"odinger  equation $ -\Delta u  =(k+V ) u$  in an
infinite slab $S= \R^{n-1}\times (0,1)$, where $V\in L^\infty(S)$ is
supported on a set $D$ of finite measure. We prove necessary
conditions for the existence of nontrivial  admissible solutions.
These conditions involve $||V||_\infty$, the measure of $D$, and the
distance   of  $k$ from the set  ${\mathcal K}= \{ \pi^2 m^2, \ m\in
\N\}$. In many cases, these inequalities are sharp.

\vspace{.5\baselineskip}


\vspace{.5\baselineskip}


\end{abstract}

\section{\bf Introduction}

\medskip
\noindent Let $S=\R^{n-1}\times (0,1)=\{(x,y)=(x_1,\cdots,x_{n-1},\,
y)\in\R^n\, :\, 0<y<1\}$ be an infinite slab.  Here, $n\ge 2$. We
study the following Schr\"odinger equation in $S$
\begin{equation}
\left\{\begin{array}{ll} -\Delta u (x,y)=(k+V(x,y))u(x,y) &\mbox{ in
}\ S\subset
\mathbb{R}^n\\
u(x)=0 &\mbox{ on }\
\partial S\end{array}\right.\label{eqk}
\end{equation}
where $k$ is a real number  and $V\in L^\infty(S)$ is supported on
some $D\subset S$ of finite measure. Unless otherwise specified, $u$
is an {\it admissible} solution in $H^1_{\rm{loc}}(S)$,  as defined
in the next section.

We are interested  in the uniqueness of solutions of \eqref{eqk}.
When $V\equiv 0$, the only admissible solution of \eqref{eqk} is the
trivial solution $u\equiv 0$ (see, for example, \cite{S, RW}). In
this paper  we prove that, unless $||V||_\infty$ is very large (in a sense to be specified)
\eqref{eqk} has only the trivial solution. The main results in this paper are
special cases of the following theorem. Together, they constitute a
proof of it.
\begin{Thm}\label{Main} Suppose that the equation \eqref{eqk} has nontrivial admissible
solutions. Then, there is a constant $c>0$,  which may depend on
$n$, $k$ and    $D$, such that
$$c ||V||_\infty \ge 1.$$
\end{Thm}

The slab is an important domain in the study of wave guides (see,
for example, \cite{BGWX}). The results for $k>0$ in
Theorem~\ref{Main} have their counterparts for acoustic wave guides
in a slab. In fact, the Helmholtz wave equation reads
\begin{equation*}
\Delta u(x,y)+ kn^2(x,y)u(x,y)=0
\end{equation*}
where $\sqrt{k}$ is the wave number, and the refraction index
$n(x,y)$ is equal to $1$ outside a compact set.  The transformation
$V(x,y)=k (n^2(x,y) -1)$ reduces the Helmholtz equation to the
Schr\"odinger equation \eqref{eqk}; the estimates
for $n(x,y)$ follow from the estimates for $V$.

Theorem~\ref{Main} gives a sufficient condition $c ||V||_\infty < 1
$ for the uniqueness of solutions of the homogeneous equation
\eqref{eqk}. Using the Lax-Phillips method (see, for example,
\cite{I})  it is proved in \cite{LU} that  this
implies the uniqueness and existence of
solutions of the nonhomogeneous equation
\begin{equation*}
\left\{\begin{array}{ll} -\Delta u (x,y)=(k+V(x,y))u(x,y)+F(x,y)
&\mbox{ in }\ S\subset
\mathbb{R}^n\\
u(x)=f(x) &\mbox{ on }\
\partial S\end{array}\right.
\end{equation*}
where $F(x,y)\in H^{-1}(S)$ has compact support in $\mathbb{R}^n$
and $f(x)\in H^{1/2}(\partial S)$ has compact support in
$\mathbb{R}^{n-1}$. So the estimate $c ||V||_\infty < 1 $ is also a
sufficient condition for the solvability of the above nonhomogeneous
equation in a slab.

Minimal potential results for the Schr\"odinger equation  $
-\Delta u= Vu$ are studied in \cite{DH, DEHL}. In these papers, the
authors consider the problems in a domain which is either bounded or
unbounded with finite measure. They show that if this equation has
nontrivial solutions, then the norm of $V$ is bounded below by a
constant. In this paper, we study the problem in an infinite slab.
We are not aware of any such uniqueness theorems for the
Schr\"odinger equation on domains of infinite measure. Because of
that, our methods are different from \cite{DH, DEHL}. Taking
advantage of the geometry of the slab, we expand the solutions with respect to the eigenfunctions    of $-\Delta-k$ in
the direction  perpendicular to the boundary
(i.e., in the $y$ direction) and then study each term. We obtain estimates that
depend on the distance between $k$ and the corresponding
eigenvalues, which are new in the literature. Also, we construct examples to show that most
of our estimates are sharp in some sense.

This paper is organized as follows.  In Section 2, we define the
admissible solutions and prove an important version of Theorem
\ref{Main} using Fourier methods, which apply in all dimensions. In
Section 3, 4 and 5 we use convolution methods to improve Theorem
\ref{Main} when $n=2$, $n=3$ and $n\ge 4$ respectively. The case
$n=3$ is probably the most interesting for applications and presents
the most technical difficulties. Section 6 is devoted to examples.

\medskip

\section{\bf  Fourier methods}
In this section  we   prove some basic estimates using Fourier
transform methods,  but first we   introduce some notation needed
throughout the paper. We will consider the admissible solutions $u$
of \eqref{eqk} which we define below using Fourier series (see also
\cite{S}). Recalling that $\{\sqrt 2\,\sin(m\pi y)\}_{m\geq 1}$
forms an orthonormal basis of $L^2(0,1) $, the solution $u$ of
\eqref{eqk} can be written in the Fourier series
\begin{equation*}
u(x,y)=\sum^\infty_{m=1}u_m(x)\sin(m\pi y)
\end{equation*}
where
\begin{equation}\label{coeff}
u_m(x)=2\int_0^1 u(x,y)\sin(m\pi y)dy.
\end{equation}
Let $f=Vu$. As done for $u$, we can write $\dsize f(x,y) =
\sum^\infty_{m=1} f_m(x) \sin(m\pi y)$.  Then by \eqref{eqk},
\begin{equation}\label{umPDE}
   -\Delta_{x} u_m+\big(m^2\pi^2-k\big)u_m=f_m
\end{equation}
 in the distribution sense, where
$\Delta_{x}=\frac{\partial^2}{\partial
x_1^2}+\cdots+\frac{\partial^2}{\partial x_{n-1}^2}$.

Assume that $V$ has support in $D=I\times (0,1) \subseteq \R^n$,
where $I\subset \R^{n-1}$ has finite measure. Note that $f_m$ has
support in $I$  because $f=Vu$ is supported in $D$.

Define
\begin{equation}\label{def-km}
k_m=\sqrt{k-m^2\pi^2}=\left\{\begin{aligned}&\sqrt{k-m^2\pi^2}\quad
& \mbox{for }m^2\pi^2< k\\
&i\sqrt{m^2\pi^2-k}\quad & \mbox{for }m^2\pi^2\ge
k.\end{aligned}\right.
\end{equation}
  So \eqref{umPDE} can be written as
\begin{equation}\label{umPDE2}
-\Delta_{x} u_m-k_m^2u_m=f_m.
\end{equation}
\begin{dfn} We say that a weak solution $u\in H_{loc}^1(S)$ of \eqref{eqk} is
an admissible solution   if
\begin{itemize}\item[1)]
  $u_m$ satisfies the radiation conditions
\begin{equation*}
u_m(x)=O(r^{\frac{2-n}{2}})\ \ \mbox{and}\ \ \left(\frac{\partial
}{\partial r}-ik_m\right)u_m(x)=o(r^{\frac{2-n}{2}})\  \mbox{as}\
r=|x|\to\infty
\end{equation*}
for $m^2\pi^2-k<0$, where $\frac{\partial u }{\partial
r}=\frac{x}{|x|}\cdot \nabla u$, and
\item[2)] $u_m\in H^1(\mathbb{R}^{n-1})$ for $m^2\pi^2-k\ge 0$.
\end{itemize}
\end{dfn}
Since the radiation conditions are imposed on $u_m$ rather than $u$,
this definition is referred to as the {\it partial radiation
condition} \cite{S}. We could consider solutions $u$ of \eqref{eqk}
in the usual Sobolev space $H_0^1(S)$, so that each $u_m\in
H^1(\R^{n-1})$. Assuming $I\subset B_{\rho}(0)=\{x: |x|<\rho\}$ and
$m^2\pi^2-k<0$, Rellich's lemma (see, for example, Lemma 35.2 in
\cite{E}) and \eqref{umPDE} imply that $u_m(x)\equiv 0$ for
$|x|>\rho$ and it satisfies the radiation condition automatically.
Thus, solutions in $H_0^1(S)$ are also admissible solutions.

Unless stated otherwise, we assume throughout this paper that $u$ is
a nontrivial admissible solution of \eqref{eqk}, and that $u_m$ is
defined as in \rq{coeff} and satisfies \rq{umPDE2}.   We denote by
$C$ a numeric constant that may change from line to line. It may
depend on the dimension $n$, but not on $D$, $u$, $m$ or
$k$.

Our strategy to prove Theorem~\ref{Main} is based on the following
simple lemma about the $u_m$.
\begin{Lemma}\label{Lmain} Suppose that there exist constants $c_m= c(k_m, D)$,
for which
\begin{equation}
||u_m||_{L^2(I)} \le c_m ||f_m||_{L^2(I)}.\label{L2norm}
\end{equation}
Then, Theorem \ref{Main} holds with  $c=\sup_m  c_m$.
\end{Lemma}
\noindent {\it Proof.}  By Parseval's formula and \eqref{L2norm},
$$||u||_{L^2(D)} \le c ||f||_{L^2(D)} =  c ||Vu||_{L^2(D)}\leq c
||V||_\infty||u||_{L^2(D)}.\ \Box
$$
\medskip
Much of the paper is devoted to finding explicit formulas  for the
$c_m$'s, which in turn provide explicit formulas for $c$. Our upper
bounds for $c_m$ often depend on the distance from $k$ to the
special set $$\mathcal{K} = \{m^2\pi^2: m = 1,2, \ldots\} .$$ Recall
that $|k_m| = \sqrt{|m^2\pi^2-k|} .$  Let
\begin{align}
&\delta_+=\min \{|k_m|: m^2\pi^2>k\}, &\delta_-=\min \{|k_m|:
m^2\pi^2<k\}.
\end{align}
When $k<\pi^2$  we let $\delta_-=\infty$.
 When $\delta_+$ and $\delta_-$
are large, Fourier methods often give the best estimates of $c_m$.
In later sections we use convolutions to estimate $c_m$ when either
$\delta_+$ or $\delta_+$ is small, but nonzero. We handle the case
$k \in \mathcal{K}$ separately.

The main result of  this section is the following.
\begin{Thm}\label{ThmX}
Assume that $k\not\in \mathcal{K}$ and  $D= I\times (0,1)$ is
bounded, with  $I \subseteq B_\rho(x_0)$ for some $\rho \delta_- >
0.1$ and $x_0\in\R^{n-1}$. Then
\begin{equation}\label{e-ThmX} \max\left\{\frac 1 {\delta_+^{2}},\ \frac{ C\rho}{\delta_-}\right\} ||V||_\infty > 1 \end{equation}
with $C$ as in Lemma \ref{Global2}.
\end{Thm}
Theorem \ref{ThmX} follows from Lemma \ref{Lmain} and  the lemmas
below.
\begin{Lemma}\label{Global1} If $\ m^2\pi^2>k$ then
\begin{equation}
\|u_m\|_{L^2(I)}\leq\frac{1}{|k_m|^2}\|f_m\|_{L^2(I)}\leq\frac{1}{\delta_+^2}\|f_m\|_{L^2(I)}.
 \label{sec2_bigm}\end{equation}
\end{Lemma}
\begin{Lemma}\label{Global2}  With $ m^2\pi^2<k$ and the assumptions of Theorem \ref{ThmX},
\begin{equation}\|u_m\|_{L^2(I)}\leq\frac{C\rho}{k_m}\|f_m\|_{L^2(I)}\leq\frac{C\rho}{\delta_-}\|f_m\|_{L^2(I)}.
\label{sec2_smallm} \end{equation}
\end{Lemma}
Example 6.4 shows
that $(\delta_+)^{-2}$ in \eqref{sec2_bigm} cannot be replaced by
$c(\delta_+)^{-2}$ for any $c<1$. When $n=2$ Example 6.7 shows that
$C (\delta_-)^{-1}$  in \eqref{sec2_smallm} cannot be improved, for
example, to $C (\delta_-)^{-\beta}$, with $\beta>1$. Likewise,
Theorem \ref{ThmX} cannot be improved in these ways.
\medskip
Lemma \ref{Global1} holds even with no assumptions on $I$, which
makes it rather unique in this paper. If $k<\pi^2$, then  Lemma
\ref{Global2} is not needed to prove Theorem \ref{ThmX}, so the
theorem holds without any assumption on the support of $V$.

\begin{Cor}\label{Csmallk}
If $k<\pi^2$ then
$$ ||V||_\infty \geq \delta_+^2 = \pi^2-k.$$
\end{Cor}

Example 6.9 shows that  $\delta_+^2$ cannot be replaced by
$\delta_+^2+\epsilon$   for any $\epsilon>0$.

\medskip
In the proof of Lemma \ref{Global1}, we let $\hat{h}=\int_{\R^{n-1}}
e^{-ix\xi} h(x)dx$ denote the Fourier transform of $h\in
L^2(\R^{n-1})$. It is well known that \eqref{umPDE2} implies
\begin{equation}\label{gmHat}
\hat{u}_m(\xi)=\frac{1}{|\xi|^2-k_m^2}\hat{f}_m(\xi) =
\hat{g}_m(\xi)\hat{f}_m(\xi)
\end{equation}
where $g_m$ is the fundamental
solution of $-\Delta_{x}-k_m^2$. So,
\begin{equation}\label{convo}
u_m(x) = g_m * f_m (x) = \int_{\mathbb{R}^{n-1}}g_m (x-y)f_m (y)dy
\end{equation}
is the convolution of $g_m$ and $f_m$. In  Sections 3, 4 and 5 we
will study $g_m$ in detail, and use \eqref{convo} and Young's
inequality for convolution to improve the estimates in  Lemmas
\ref{Global1}  and \ref{Global2} for small values of $\delta_-$ and
$\delta_+$.
\medskip
\noindent {\it Proof of Lemma~\ref{Global1}}.  To prove
\eqref{sec2_bigm}, we observe that $k_m^2 = -|k_m|^2$, and by
\eqref{gmHat},
$\|\hat{u}_m\|_{L^2(\mathbb{R}^{n-1})}\leq\frac{1}{|k_m|^2}\|\hat{f}_m\|_{L^2(\mathbb{R}^{n-1})}$.
By Plancherel's theorem, and because    $f_m$ has support  in $I$,
\begin{equation*}
\|u_m\|_{L^2(I)}\leq\|u_m\|_{L^2(\mathbb{R}^{n-1})}\leq\frac{1}{|k_m|^2}\|f_m\|_{L^2(\mathbb{R}^{n-1})}
\leq\frac{1}{\delta_+^2}\|f_m\|_{L^2(I)}.\ \Box
\end{equation*}
\medskip
\noindent {\it Proof of Lemma \ref{Global2}}.  To prove
\eqref{sec2_smallm}, we can assume $x_0=0$ without loss of
generality   and $\rho=1$ by dilation (see the remarks following
\rq{diln}). Note that $k_m^2>0$ and the term
$\frac{1}{|\xi|^2-k_m^2}$ in \eqref{gmHat} is not integrable. We use
\begin{equation*}
\frac{1}{|\xi|^2-k_m^2-i0}=\lim_{\epsilon\to
0^+}\frac{1}{|\xi|^2-k_m^2-i\epsilon}
\end{equation*}
as a regularization of $\frac{1}{|\xi|^2-k_m^2}$. Since $u_m$
satisfies the radiation condition, it can be represented as
\begin{equation*}
u_m(x)=\frac{1}{(2\pi)^{n-1}}\int_{\mathbb{R}^{n-1}}\frac{e^{ix\cdot\xi}\hat{f}_m(\xi)}{|\xi|^2-k_m^2-i0}d\xi.
\end{equation*}
(See, for example, Theorem 19.5 in \cite{E}). Denote by
$W_{0,s}(\mathbb{R}^{n-1})$ the weighted Sobolev space with the norm
\begin{equation*}
\|g\|_{0,s}=\left(\int_{\mathbb{R}^{n-1}}(1+|x|)^{2s}|g(x)|^2dx\right)^{\frac{1}{2}}.
\end{equation*}
Since $f_m$ has support in $I\subset \mathbb{R}^{n-1}$, we know
$f_m\in W_{0,s}(\mathbb{R}^{n-1})$ for any $s$. Using the Agmon's
estimate (see, for example, Theorem 29.1 in \cite{E}) we have
\begin{equation*}
\|u_m\|_{0,-1}\le \frac{C}{k_m}\|f_m\|_{0,1}
\end{equation*}
 for $k_m\geq 0.1$. Since $f_m$ has support in $I\subset B_1(0)$, we
have
\begin{eqnarray*}
\|u_m\|_{L^2(I)}&=&\left(\int_{I}|u_m|^2dx\right)^{\frac{1}{2}}
\leq\left(\int_{I}\left(\frac{1+|x|}{2}\right)^{-2}|u_m|^2dx\right)^{\frac{1}{2}} \\
&=&2 \|u_m\|_{0,-1}\leq
2\frac{C}{k_m}\|f_m\|_{0,1}=2\frac{C}{k_m}\left(\int_{I}(1+|x|)^{2}|f_m|^2dx\right)^{\frac{1}{2}}\\
&\leq& 2\frac{C}{k_m}\cdot
2\left(\int_{I}|f_m|^2dx\right)^{\frac{1}{2}} =
\frac{4C}{k_m}\|f_m\|_{L^2(I)}\leq\frac{4C}{\delta_-}\|f_m\|_{L^2(I)}.
\end{eqnarray*}
Now $4C$ serves as the generic constant in part (2). This concludes
the proof of the lemma. $\hfill\Box$

\section{\bf Results in dimension 2 and some general lemmas}

We now present versions of Theorem \ref{Main} in dimension $n=2$
which do not depend on Fourier methods, but rather on convolution.
This section also includes some remarks on dilation and
rearrangement which apply in every dimension.

\subsection{The main results in dimension 2}

We now estimate the $L^2$ norm of $u_m  = g_m*f_m$  (see
\eqref{convo}) using Young's inequality for convolution. We present
the case $n=2$ first, to show the basic ideas of this approach most
clearly. The formula for $g_m$ is relatively simple in this case.
When $k_m=0$,  \rq{umPDE} takes the form $-u_m''=f_m$, and its
fundamental solution is $g_m(x)= -\frac{|x|}2$. Otherwise, with
$|x|=r$, the equation \eqref{gmHat} implies (see, for example,
\cite{RW})
\begin{equation}\label{egm_n2}
g_m(x)=g_m(r) =-\frac{e^{ik_mr}}{2ik_m}.
\end{equation}
We use this
to improve on Lemmas \ref{Global1} and \ref{Global2} and Theorem
\ref{ThmX} when $\delta =\min\{\delta_+,\ \delta_-\}$ is small.

\begin{Thm}\label{n2main}   Suppose $n=2$.

1) If $k\not\in \mathcal{K}$, then \BEL{n2main1}{||V||_\infty
|D|\over 2\delta}  \geq 1.\EE

2) If $k\in \mathcal{K}$ and $D\subseteq B_\rho(0)\times (0,1)$ with
$\rho \ge 1$, then \BEL{n2main2} 2\rho ||V||_\infty |D|  \geq 1.\EE
\end{Thm}

Examples 6.5 and 6.11 show that these inequalities are sharp when $\delta$ is small.
But if $\delta$ is large enough, Theorem \ref{ThmX} may give a
better result than Theorem \ref{n2main}. To see this, suppose that
$\delta_- |D| > {2C\rho\delta_+ } > \delta_+ |D| >2$, so that
$\delta =\delta_+$. Also, ${|D|\over 2\delta} > \delta_+^{-2}$ and
${|D|\over 2\delta} > {C\rho\over \delta_-}$, so in this case
\eqref{e-ThmX}  is stronger than \rq{n2main1}.
\medskip
\noindent {\it Proof:} Part 1) follows directly from Lemmas
\ref{Lmain} and \ref{n2lemma}. For part 2), note that if
$m^2\pi^2=k$, the $c_m$ defined by Lemma \ref{n2Spec} is $2\rho
|D|$. For other $m$, our estimate of $c_m$ from Lemma \ref{n2lemma}
is $\frac{|D|}{2|k_m|}$, which is smaller, since $\rho$, $|k_m|\ge
1$. Part 2) now follows from Lemma \ref{Lmain}. $\hfill\Box$

\begin{Lemma}\label{n2lemma}    Suppose $n=2$. If $k_m
\not=0$ then
\begin{equation}\label{est_um_n2}
||u_m||_{L^2(I )} \leq \frac{|I|}{2|k_m|}||f_m||_{L^2(I )} \leq
\frac{|I|}{2\delta} ||f_m||_{L^2(I )}.
\end{equation}
\end{Lemma}
\medskip
\noindent {\it Proof:} By \rq{egm_n2}, $|g_m(x)|\leq
\frac{1}{2|k_m|}$. So, $||u_m||_{L^\infty(I)} \leq  \frac{1}{2|k_m|}
||f_m||_{L^1(I)}$. Since $||u_m||_{L^2(I)} \le |I|^{\frac
12}||u_m||_{L^\infty(I)}$, $||f_m||_{L^1(I)} \le |I|^{\frac
12}||f_m||_{L^2(I)}$ and $\delta \le |k_m|$, the result follows.
$\hfill\Box$

Examples \ref{2Ex} and \ref{n2mD} show that the inequality
\eqref{est_um_n2} is sharp.

\begin{Lemma}\label{n2Spec}   With the assumptions of Theorem \ref{n2main} part (2), and with $m^2\pi^2=k$,
\BEL{n2Sp} \|u_m\|_{L^2(I)}\leq 2\rho |I|\ \|f_m\|_{L^2(I)}. \EE
\end{Lemma}
\noindent {\it Proof:} As mentioned above, since $k_m=0$, we have
$u_m = f_m*g_m$, with $g_m(x)= -\frac{|x|}2$, but the fundamental
solution is not unique in this case. Note that $\int_{\R} f_m\, dx =
-\int_\R u_m''\,dx =0$.  So, $u_m = f_m*g_m = f_m*(\rho+g_m)$. The
function $G_m = \rho -\frac{|x|}2 $ is radially decreasing and
nonnegative on $(-2\rho,\ 2\rho)$ which contains the set $2B$ (a
ball of radius $|I|$ centered at 0). Lemma \ref{rearr} gives
$$\|u_m\|_{L^2(I)}\leq \|G_m \|_{L^1(2B)} \|f_m\|_{L^2(I)}$$
and $\|G_m \|_{L^1(2B)} \le 2\rho |I|$. $\hfill\Box$

Example \ref{n2Sp*} shows the inequality \eqref{n2Sp} is sharp.

\medskip
As stated in Lemma~\ref{Lmain}, our results are mainly based on the
estimates of the form $||u_m||_{L^2(I)} \le c_m ||f_m||_{L^2(I)}$.
If we have two such estimates, with different bounds $c_m$ and
$c_m'$, we can use the minimum of the two. We can combine Theorems
\ref{ThmX} and \ref{n2main} in this way for a better result.

 \begin{Thm}\label{T-combined} Suppose $n=2$. If $k\not\in \mathcal{K}$, then
 \begin{equation}\label{E-MaxV2}
c||V||_\infty\geq 1, \quad \mbox{where}\ c=
\max\left\{\frac{|D|}{2\delta_-}, \
\min\left\{\frac{|D|}{2\delta_+},\frac{1}{\delta_+^2}\right\}\right\}.
\end{equation}
 \end{Thm}

\noindent {\it Proof.} When $m^2\pi^2>k$, we have the estimate
\eqref{sec2_bigm} in Lemma~\ref{Global1} and \eqref{est_um_n2} in
Lemma~\ref{n2lemma}. Thus,
$$
||u_m||_{L^2(I )} \leq
\min\left\{\frac{|I|}{2|k_m|},\frac{1}{k_m^2}\right\} ||f_m||_{L^2(I
)}\leq \min\left\{\frac{|I|}{2\delta_+},\frac{1}{\delta_+^2}\right\}
||f_m||_{L^2(I )}.
$$
When $m^2\pi^2<k$, we use the estimate \eqref{est_um_n2}, and also
$$
||u_m||_{L^2(I )} \leq \frac{|I|}{2|k_m| } ||f_m||_{L^2(I )}\leq
\frac{|I|}{2\delta_- } ||f_m||_{L^2(I )}.
$$
 These estimates  and  the fact that  $|D|=|I|$ yield
\eqref{E-MaxV2}. $\hfill\Box$

\medskip
 \noindent
{\it Remark.} In the case $\frac{|D|}{2\delta_-}\leq
\frac{1}{\delta_+^2} <\frac{|D|}{2\delta_+}$, we have
$c=\frac{1}{\delta_+^2}$ in \eqref{E-MaxV2}. The constant in
\eqref{n2main1} is $\frac{|D|}{2\delta_+} >c$, so \eqref{E-MaxV2} is
stronger than  \eqref{n2main1}. In all other cases, the constants
are the same. The results in other dimensions can be written in the
form of Theorem \ref{T-combined}, but they tend to be even more
complex, so we leave those forms to the interested reader.

\subsection{Dilation and rearrangement.}

We describe when an inequality on $||u_m||_2$ of the form
\rq{L2norm} is dilation-invariant. Suppose \BEL{e-diff-dilat}-\Delta
w(x)  - h^2 w(x) = g(x),\EE with $g$ supported on $I\subset
\R^{n-1}$. Let $W(x) = w(\lambda x)$ and $G(x) = g(\lambda x)$. Then
\BEL{diln}-\Delta W(x)  - (\lambda h)^2 W(x) = \lambda^2 G(x)\EE
with $G$ supported on $J = \lambda^{-1}I$.

Suppose, for example, that Lemma \ref{Global2}  holds when $\rho=1$.
Suppose $w = u_m$ satisfies \rq{umPDE} (which is
\eqref{e-diff-dilat} with $h=k_m$, and $g=f_m$), but that $f_m$ is
supported in $I=B_\rho(0)$ for some $\rho\not=1$.
To prove that Lemma \ref{Global2} also holds for $w$, we
define $W$ as above with $\lambda = \rho$. Since $G$ is supported on $B_1(0)$,
and \eqref{diln} holds, we can apply Lemma \ref{Global2} to $W$, to obtain
$$||W||_{L^2(J)}\leq \frac{C}{\rho k_m}||\rho^2 G||_{L^2(J)},$$
and then
$${||w||_{L^2(I)} \over ||g||_{L^2(I)}} = {||W||_{L^2(J)} \over ||G||_{L^2(J)}} \le {C \over \rho k_m} \rho^2,$$
which shows the lemma holds for $w=u_m$. To summarize, we may assume
$\rho =1$ in the proof of Lemma \ref{Global2} and in the sharpness
Example 6.7. Similar remarks apply to all our lemmas of this type,
where the ratio of the norms of $u_m$ and $f_m$ has an upper bound
of the form $\rho^2 \phi(\rho k_m)$ for some function $\phi$.

In other lemmas, $I$  is assumed to have finite measure, but is not
necessarily bounded. By setting $\lambda = |I|^{1\over n-1}$ and
applying similar reasoning, we get $|{\rm supp\ }G|=1$ and can study
$||W||_{L^2(J)}/||G||_{L^2(J)}$ as done above. Such lemmas are
dilation-invariant if $||u_m||_{L^2(I)}/||f_m||_{L^(I)}$ has an
upper bound of the form $|I|^{2\over n-1} \phi(|I|^{1\over n-1}
k_m)$. For example, in  Lemma \ref{n2lemma} this occurs with $\phi
(t) = {1\over 2|t|}$.

When $k_m=0$, $\phi$ is constant. For example, note that the ratio
in Lemma \ref{n2Spec} is neither $\rho^2$ not $|I|^2$, but it still
has the correct homogeneity for dilation invariance. Note that we
cannot apply these dilation remarks to results in the form of
Theorem \ref{Main}.

\medskip
We now prove a rearrangement lemma. Recall that $u_m=f_m*g_m$, (see
\eqref{convo}) and that $f_m$ is supported in $I$.

\begin{Lemma}\label{rearr}
Suppose $I^{*}=I-I \subseteq B_R(0) \subseteq \R^{n-1}$, with
$0<R\le \infty$. Suppose  $|g_m| \le G_m = \overline{G_m}$ on
$B_R(0)$, where $\overline{G_m}$ is the symmetric non-increasing
rearrangement of $G_m$. Then $$||u_m||_{L^2(I)} \leq
\|G_m\|_{L^1(2B)} ||f_m||_{L^2(I)}$$   where $B$ is the ball in
$\mathbb{R}^{n-1}$ centered at $0$ with the same measure as $I$.
\end{Lemma}

In some applications we do not have $G_m = \overline{G_m}$ on
$B_R(0)$, but can get it easily by redefining $G_m$ to be zero
outside $B_R(0)$ before applying this lemma.
\medskip
\noindent {\it Proof.} Recall the well known rearrangement
inequality by Hardy and Littlewood (see, for example, \cite{LL} pg.
76). If $f$, $g$ and $h$ are non-negative then
\begin{equation}\label{EHL}
\int_{\R^{n-1}} h(x) [f*g(x)]\ dx \le \int_{\R^{n-1}} \bar{h}(x)
[\bar{f}*\bar{g}(x)]\ dx.
\end{equation}
Clearly $||\bar{h}||_{L^p(B)} = ||h||_{L^p(I)}$ for all $p$. We
apply \eqref{EHL} with $h = \chi_I| u_m|$, $g=G_m$ and $f=|f_m|$.
With Holder's inequality and Young's inequality we have
\begin{eqnarray*}
&& ||u_m||_{L^2(I)}^2 = \int_{\R^{n-1}}\chi_I u_m(x) (f_m*g_m)(x)
dx\\ &\leq& \int_{\R^{n-1}}\chi_I |u_m(x)| (|f_m|*|g_m|)(x)| dx\\
&\leq&\int_{\R^{n-1}}\chi_I |u_m(x)| (|f_m|*G_m)(x)| dx
\leq\int_{B}\overline{|u_m(x)|} (\overline{|f_m|}*\overline{G_m})(x)
dx \\ &\leq&
\|\overline{|u_m|}\|_{L^2(B)}\|\overline{|f_m|}\|_{L^2(B)}\|\overline{G_m}\|_{L^1(2B)}
\\ &= &||u_m||_{L^2(I)}||f_m||_{L^2(I)}\|G_m\|_{L^1(2B)}. \ \Box
\end{eqnarray*}

\section{\bf Results  in dimension $ 3$}

In this section we discuss uniqueness results in dimension $n=3$. We
provide technical estimates on the fundamental solutions $g_m$ in
Section 4.1, and complete the proofs of all lemmas in Section 4.2.

As stated in Section~2, the results in Theorem~\ref{ThmX} hold  in
every dimension, but the estimates when $m<\sqrt k/\pi$  are valid
only when $D$ is bounded. When $D$ is not bounded, we replace  the
estimates  in  Lemma~\ref{Global2} with those in Lemma \ref{Lorenz3}
below. Inequality \rq{eq_Lorenz3} also holds when $m>\sqrt k/\pi$.
In what follows, we assume that $n=3$, that \eqref{eqk} has an
admissible nontrivial solution $u$ and that $|I| =|D|$ is finite.

\begin{Lemma}\label{Lorenz3} If $k_m\neq 0$, then
\BEL{eq_Lorenz3} ||u_m||_{L^2(I)} \leq C  |k_m|^{-\frac
12}|I|^{\frac 3 4} ||f_m||_{L^2(I )}.\EE
\end{Lemma}

Lemma \ref{Lorenz3} and Lemma~\ref{Lmain} give:
\begin{Thm}\label{n3thm1}If $k\not\in \mathcal{K}$, then
\BEL{eq_n3thm1}\frac{C||V||_\infty |D|^{{\frac 3 4}}}{\delta^{\frac
12}}  \geq 1. \EE
\end{Thm}

Theorem \ref{n3thm1}  is a suitable replacement for Theorem
\ref{ThmX} when neither $|D|$ nor  $\delta$ is too large or too
small. However, we do not claim that Lemma \ref{Lorenz3} or Theorem
\ref{n3thm1} is sharp. When either  $|D|$ or $\delta$ is fixed and
the other is sufficiently large or small, Theorem  \ref{n3thm1} can
be improved. We demonstrate this in all four cases.
\medskip

\noindent{\it Case 1:  Suppose $|D|$ is fixed and $\delta$ is
sufficiently large.} When $m<\sqrt k/\pi$, we use
Lemma~\ref{Lorenz3} to estimate $c_m$. If $m>\sqrt k/\pi$, we use
Lemma~\ref{Global1}, which is stronger for large enough
$\delta$. By Lemma \ref{Lmain}, \BEL{eq_n3thm2}
\max\left\{\frac{1}{\delta_+^{2}}, \frac{C |D|^{{\frac 3
4}}}{\delta_-^{\frac 12}}\right\} ||V||_\infty \geq 1. \EE
\medskip

\noindent{\it Case 2: Suppose $|I|=|D|$ is fixed and $\delta$ is
sufficiently small.} Then,  Lemma~\ref{Global1}, Lemma~\ref{Global2}
and Lemma~\ref{Lorenz3} are not sharp, but the following estimate is
(see Example \ref{Forn3S}).
\begin{Lemma}\label{n3S}
Assume that  $0< 4 \pi^{-\frac 12} |k_m|  |I|^{\frac 12}  < 1 $.
Then \BEL{eq_L-n3S} ||u_m||_{L^2(I )} \leq C|I|(1-\ln(|I|^{\frac
12}|k_m|))||f_m||_{L^2(I )}. \EE
\end{Lemma}

For the terms $u_m$ with $0< 4 \pi^{-\frac 12} |k_m|  |I|^{\frac 12}
< 1 $, we use Lemma~\ref{n3S} to estimate $c_m\leq
C|D|(1-\ln(|D|^{\frac 12}\delta))$. For other terms, in view of $4
\pi^{-\frac 12} |k_m|  |I|^{\frac 12} > 1 $, Lemma ~\ref{Lorenz3}
shows $c_m \le C|D|^{\frac 34} |k_m|^{-\frac 12} \le C|D|$,
which is smaller than $C|D|(1-\ln(|D|^{\frac
12}\delta))$ when $\delta$ is small. By Lemma \ref{Lmain},
\BEL{n3ln} C|D|(1-\ln(|D|^{\frac 12}|\delta|)) ||V||_\infty \geq 1.
\EE

\medskip
\noindent{\it Case 3: Suppose $k$ is fixed and $|D|$ is sufficiently
small.} In view of $C|D|\leq C|D|(1-\ln(|D|^{\frac 12}\delta))$ for
small $|D|$, the same reasoning as in
Case 2 leads again to \rq{n3ln}.
\medskip
\noindent{\it Case 4: Suppose $k$ is fixed and $|D|$ is sufficiently large.}
 If $m>\sqrt k/\pi$, Lemma~\ref{Global1} is stronger than
Lemma~\ref{Lorenz3} for large $|D|$. Then the same reasoning as in
Case 1 leads to \rq{eq_n3thm2}.

\begin{Lemma}\label{n3Spec}Suppose that $k=m^2\pi^2$ and $I\subset B_\rho(0)$. Then
\BEL{eq_n3spec}||u_{m}||_{L^2(I )} \leq
\frac{|I|}{\pi}(1+\ln(\pi\rho^2/|I|))||f_{m}||_{L^2(I )}.\EE
\end{Lemma}

This lemma is sharp (see Example \ref{n3Sp*}) and allows us to
handle a final case, that $k\in \mathcal{K}$.
\begin{Thm}\label{n3thmspec}
If $k\in \mathcal{K}$ and $D\subseteq B_\rho(0)\times (0,1)$, then
\BEL{eq_n3thmspec} \max\left\{C|D|^{{\frac 3 4}},
\frac{|D|}{\pi}(1+\ln(\pi\rho^2/|D|))\right \} ||V||_\infty \geq 1.
\EE
\end{Thm}

\medskip
\noindent {\it Proof.}  If $m^2\pi^2=k$, use \eqref{eq_n3spec} to
get $c_m$. For other $m$,  the $c_m$ obtained from
Lemma~\ref{Lorenz3} are at most $C|D|^{{\frac 3 4}}$.   Lemma
\ref{Lmain} concludes the proof.

$\hfill\Box$

\medskip
\noindent{\it Remark: } If  $k\in {\mathcal K}$ and $m\ne \sqrt
k/\pi$, then $|k_m| \ge \sqrt 3\pi$. We can combine
Lemma~\ref{n3Spec} with Lemmas~\ref{Global1} and \ref{Global2}  to
obtain $$\max\left\{\frac{1}{3\pi^2}, \frac{C\rho}{\sqrt{3\pi^2}},
\frac{1}{\pi}|D|(1+\ln(\pi\rho^2/|D|)) \right\}||V||_\infty \geq
1.$$

\medskip
In the next subsections we prove Lemmas \ref{Lorenz3}, ~\ref{n3S}
and ~\ref{n3Spec} along with some estimates needed in Section 5,
where $n\ge 4$. We now assume $n\ge 3$, unless stated otherwise.

\subsection{Properties of fundamental solutions}
Recall from Section~2 that we can estimate the constants $c_m$ in
Lemma \ref{Lmain} using  the formula $u_m=f_m*g_m$, where $g_m$ is
the fundamental solution of $-\Delta_{x}-k_m^2$. By Young's
inequality for convolutions
$$
||u_m||_{L^2(I)}= ||f_m*g_m||_{L^2(I)} \leq ||g_m||_{L^1(I^*)}
||f_m||_{L^2(I)}
$$
where $I^*= I-I$. In this section we give explicit estimates for
$||g_m||_{L^1(\R^{n-1})}$.

It is well known (see \cite{Er}) that
\begin{equation}\label{Egm_sec3}
g_m(r) = {i\over 4 }\left({k_m \over 2\pi r}\right)^s H_s( k_m r) =
{i\over 4 }\left({k_m^2 \over 2\pi }\right)^s ( k_m r)^{-s}H_s( k_m
r)
\end{equation}
where $r=|x|$, $s=\frac{n-3}{2} \ge - \frac 12$, and $ H_s$ is a
Hankel function. Note $ H_s =  J_s +iY_s$ where $J_s$ and $Y_s$ are
Bessel functions of the first and second kind, resp. We have a
Poisson type representation formula for $g_m$.

\begin{Lemma}\label{poissonRep}
For $s\geq 0$, define
\begin{equation}\label{EIs}
I_s (z) = \int_0^\infty e^{-t} t^{s-\frac 12}\left(1
+\frac{t}{2z}\right)^{s-\frac 12} dt.
\end{equation}
Then
\begin{equation}\label{Eintgm}
g_m(r)=c(s)(-ik_m)^{s-\frac{1}{2}}r^{-s-\frac{1}{2}}e^{ik_m r}
I_s(-ik_m r)
\end{equation}
where $c(s)= \frac{1}{4  (2\pi)^s} \frac{ 1}{ \Gamma(s+{\textstyle
\frac 12}) } \left(\frac 2{\pi }\right)^{\frac 12}$.
\end{Lemma}

\noindent{\it Proof.}   The following representation formula is well
known (see e.g. \cite{Er}, vol. 2)
\begin{equation}\label{e-intHs}
H_s (iz)= \frac{ e^{-\frac{i\pi s}2}}{i\Gamma(s+{\textstyle \frac
12}) }  \left(\frac 2{\pi z}\right)^{\frac 12} e^{-z} I_s(z).
\end{equation}
So, by \eqref{Egm_sec3},
\begin{eqnarray*}
g_m(r) &=&{i\over 4 }\left({k_m \over 2\pi r}\right)^s \frac{
e^{-\frac{i\pi s}2}}{i\Gamma(s+\frac{1}{2})}  \left(\frac 2{\pi
(-ik_m r)}\right)^{\frac 12} e^{ik_m r} I_s(-ik_m r)\\
&=&c(s)(-ik_m)^{s-\frac{1}{2}}r^{-s-\frac{1}{2}}e^{ik_m r} I_s(-ik_m
r). \  \Box
\end{eqnarray*}

\begin{Lemma}\label{HankS}
Let $n=3$ and $m\ne \sqrt k/\pi$.

\noindent 1) If $2|k_m|r>1$,
\begin{equation}\label{e-largeR}
|g_m(r)| \leq C(|k_m|r)^{-\frac{1}{2}}.
\end{equation}
2) If $2|k_m|r \leq 1,$
\begin{equation}\label{e-gmS}
| g_m(r)| \leq  C(1 - \log (2r|k_m|)).
\end{equation}
\end{Lemma}

\noindent {\it Proof.}  By \eqref{Eintgm} with $s=0$, $g_m(r)=c(0)
e^{ik_m r} I(-ik_m r)$, where
\begin{equation*}I(z) = z^{-\frac 12} I_0(z) =
\int_0^\infty e^{-t} t^{ -\frac 12}\left( z+ \frac t {2 } \right)^{
-\frac 12} dt .\end{equation*} Either $z=-ik_mr>0$ or $k_m>0$,
implying $|e^{-z}|\le 1$, and
$$|g_m(r)| \leq c(0) |I_0(z)|\le c(0) \int_0^\infty e^{-t} t^{ -\frac 12}\left|z+ \frac t2\right|^{ -\frac 12} dt= c(0) J(z).$$
Let $J  = J_A+J_B$ where $J_A = \int_0^{2|z|} e^{-t} t^{ -\frac
12}|z+ \frac t2|^{ -\frac 12} dt$. In $J_A$, we use $|\frac t {2}
+z| \ge |z|$ to get $J_A \le |z|^{-\frac 12} \int_0^{2|z|} e^{-t}
t^{ -\frac 12}\ dt$. In $J_B$, we use $|\frac t {2} +z| \ge
\frac{t}{2}$ to get $J_B \le \sqrt 2 \int_{2|z|}^\infty e^{-t}
t^{-1} \ dt$.

Suppose ${2|z|} \ge 1$. Then $ J_A < |z|^{-\frac 12}\int_0^{\infty}
e^{-t} t^{ -\frac 12}dt= |z|^{-\frac 12}\sqrt \pi $, and $J_B \le
\sqrt 2 \int_{2|z|}^\infty e^{-t} \ dt \leq \sqrt 2 e^{-2|z|} \leq
e^{-1} |z|^{-\frac 12}$ because the function $t\to te^{1-t^2}$ is
decreasing when $t\ge 1$. So, $J_A+J_B \le 5 |z|^{-\frac 12}$ and
\eqref{e-largeR} follows.

Now, suppose $2|z| < 1$. Then $ J_A \leq|z|^{-\frac 12}\int_0^{2|z|}
t^{ -\frac 12} dt  <4.$ Note that $\int_{2 |z|}^1 e^{-t} t^{-1} \ dt
\le \int_{2|z|}^1 t^{-1} \ dt = -\ln({2|z|})$ so $J_B \le C( 1
-\ln({2|z|}))$, giving \rq{e-gmS}. $\hfill\Box$

\medskip
In Section 6, we will construct sharpness examples based on $g_m$
and the following Lemma.

\begin{Lemma}\label{Lg1/g}  For every $m \ne \sqrt k/\pi$,  and every
$a>0$,
\begin{itemize}\item [1)]
$\displaystyle |g'_m(a)|/|g_m(a)|= |k_m|$ when $n=2$, and

\item [2)]
$\displaystyle |g'_m(a)|/|g_m(a)| \le C (a^{-1}+|k_m| )$
  when $ n\geq 3$.
\end{itemize}
\end{Lemma}

\noindent {\it Proof.} When $n=2$, $g_m(x)= (2i k_m)^{-1} e^{i k_m
x}, $  so   $|g'_m(a)|/|g_m(a)|= |k_m|$.

When $n\geq 3$, we let $s=\frac{n-3}{2}\geq 0 $.  As always,  $C$
will denote a generic constant that may change from line to line. We
assume $m>\sqrt k/\pi$  (so $-ik_m= |k_m|$) since the proof is
similar in the other case.

Recalling that   $\frac{d}{dz} (z^{-\alpha} H_\alpha(z))= -
z^{-\alpha} H_{\alpha+1}(z)$ (see e.g. \cite{Er}),   and the formula
for  $g_m$ in \eqref{Eintgm},
\begin{eqnarray*}
&&g_m'(r) =  C |k_m|^{2s} \frac{d}{dr} ((ir|k_m|)^{-s} H_s(ir|k_m|))\\
&=&C |k_m|^{2s }  ( i|k_m|)\frac{d}{dz}(z^{-s} H_s(z)) =
C|k_m|^{2s+1}z^{-s} H_{s+1}(z).
\end{eqnarray*}
Recalling the estimates of the  Hankel functions
$$|H_\alpha(z)|\leq C(\alpha)|e^{-z}| |z|^{-\alpha}(1+
|z|^{\alpha-\frac 12}), $$ which are valid for $\alpha>0$ and
$Arg(z) <\pi$, we can see that
\begin{eqnarray}\label{e-g'm(1)}
&&|g_m'(a)|= C   |k_m|^{2s+1} (a|k_m|)^{-s}
|H_{s+1}(i|k_m|a)|\nonumber\\
&\leq&   C |k_m|^{ 2s+1} ( a|k_m|)^{-s} e^{-|k_m| a} (a
|k_m|)^{-s-1} (1+ (a|k_m|)^{s+\frac 12}) \nonumber\\
&=& C e^{-|k_m| a}a^{-2s-1}(1+ (a|k_m|)^{s+\frac 12}).
\end{eqnarray}
We now find a lower bound for $|g_m(a)|$. From the integral formula
\eqref{Egm_sec3}, follows that
 \begin{equation}\label{e-intgm}
|g_m(a)| = C e^{-|k_m| a } a^{-2s  } \int_0^\infty\!\!\! e^{-t}
t^{s-\frac 12}\left(|k_m| a  +\frac{t}{2  }\right)^{s-\frac 12} dt.
\end{equation}
We assume  $n\ge 4$. The proof   for $n=3$  is very similar, and it
is left to the reader.   When $a|k_m|\leq 1$,
$$
|g_m(a)| \ge   C e^{-|k_m| a }   a^{-2s  } \int_0^\infty\!\!\!
e^{-t} t^{ 2s-1}  dt = C e^{- a|k_m| } a^{-2s  }
$$
and so, from \rq{e-g'm(1)},
\begin{equation}\label{e-small}
|g'_m(a)|/|g_m(a)|  \leq  C \frac{  e^{-|k_m| a }a^{-2s-1} (1+
(a|k_m|)^{s+\frac 12})}{e^{-|k_m| a } a^{-2s  }}\leq  \frac{C}{a}.
\end{equation}
When $a|k_m|\ge 1$,
$$
|g_m(a)| \ge \!  C e^{-|k_m| a }  a^{-2s} (a|k_m|)^{s-\frac 12}
\int_0^\infty\!\!\! e^{-t} t^{s-\frac 12} dt= C e^{-|k_m|a }
|k_m|^{s-\frac 12} a^{-s-\frac 12}
$$
 from which follows that
\begin{align}\label{e-large}
& |g'_m(a)|/|g_m(a)|  \leq  C \frac{  e^{-|k_m| a} a^{-2s-1} (1+
(a|k_m|)^{s+\frac 12})}{e^{-|k_m| a } |k_m|^{s-\frac 12} a^{-s-\frac
12} } \nonumber\\
&= C \frac{    (1+ (a|k_m|)^{s+\frac 12})}{ |k_m|^{s-\frac 12}
a^{s+\frac 12} } =  C|k_m| ((a|k_m|)^{-\frac 12-s}+ 1)\leq C
|k_m|. \
\end{align}
So, when $n\ge 4$  we have proved that
    $|g'_m(a)|/|g_m(a)| \le C\max\{  |k_m|, \frac 1a\} \leq  C(a^{-1}+ |k_m|),
$ as required.  $\hfill\Box$

\subsection{Proofs of the main lemmas} In this subsection we
complete the proofs of Lemma~\ref{Lorenz3}, Lemma~\ref{n3S} and
Lemma~\ref{n3Spec}.
\medskip
\noindent{\it Proof of Lemma~\ref{Lorenz3}.} To prove this result,
we need the following version of Young's inequality
 in  Lorentz spaces (See e.g. \cite{G})
for  the properties of   these spaces)
\begin{equation}\label{e-Young} ||u_m||_{4,\infty}\le C ||f_m||_{1,\infty}||g_m||_{4,\infty}.
\end{equation}
The inequality  \eqref{e-Young}  and
 \BEL{c}
   ||g_m||_{4, \infty}\leq C |k_m|^{-\frac 12},\EE
\BEL{a} ||f_m||_{1,\infty} \le  C ||f||_{1} \leq C |I|^{\frac
12}||f_m||_2 ,  \EE \BEL{b} ||u_m||_2 \le    C|I|^{{\frac 14}}
||u_m||_{4,\infty}    \EE imply \eqref {eq_Lorenz3}.

We prove  \eqref{c} and we leave the proof of   \eqref{a} and
\eqref{b}   to the reader. Since $|1-\ln(t)| \le C t^{-1/2}$ for
$0<t<1$, Lemma \ref{HankS} shows $|g_m (r)| \leq C |k_m r|^{-1/2}$
for all $r$. Note $|\{ x\in \R^2: |k_m r|^{-1/2} > \alpha \}|  =
\alpha^{-4}|k_m|^{-2}$ and $\alpha(\alpha^{-4}|k_m|^{-2})^{-1/4} =
|k_m|^{-\frac 12}$, proving \eqref{c}. $\hfill\Box$

\medskip

\noindent{\it Proof of Lemma~\ref{n3S}.} We first find a suitable
$G_m$ in order to use rearrangement Lemma \ref{rearr}. By Lemma
\ref{HankS}, the following estimate is valid for every $r>0$ and
every $m\ne \sqrt k/\pi$ for which $2r|k_m|<1$,
$$|g_m(r)| \leq G_m(r) = C(1 - \ln (2r|k_m|)).
$$
By Lemma \ref{rearr}, it is enough to show that  $||G_m||_{L^1(2B)}
\le C|I|(1+|\ln(|I|^{ \frac 12} |k_m| )|)$ where $|2B|=4|I|$. This
implies the radius of $2B$ is $R$, with $\pi R^2 = 4|I| $ and by
assumptions
$$  4 |k_m|^2 R^2    = 16 \pi^{-1} |k_m|^2 |I|  < 1.  $$
 Thus, $2r|k_m| <2R|k_m|<1$, and
\begin{eqnarray*}||G_m||_{L^1(2B)} &\le&  C \int_0^{R} r(1-\ln(2r|k_m|) ) \ dr\\
&\le & C   R^2 (3-2 \log (2|k_m| R)) \le C |I|(1+ |\ln(|I|^{\frac
12} |k_m| )|)\end{eqnarray*} as desired.$\hfill\Box$

\medskip

\noindent{\it Proof of Lemma~\ref{n3Spec}.} Since $V$ has compact
support, $-\Delta u_m = f_m$ has compact support. By Rellich's
lemma, $u_m$ also has compact support. By the divergence theorem,
$\int_{\mathbb{R}^2}f_mdx=-\int_{\mathbb{R}^2}\Delta u_mdx=0$. So,
for any constant $c$, we have $c*f_m\equiv 0$. Let
$G_m=\frac{1}{2\pi}\ln(2\rho/|x|)$. So $u_m=G_m*f_m$ and $G_m$ is a
radially decreasing positive function on $B_{2\rho}(0)$, which
contains the set $2B$, with $B$ as in Lemma \ref{rearr}. It is
enough to show that
\begin{equation*}
||G_m||_{L^1(2B)} \le \frac{1}{\pi}|I|(1+\ln(\pi\rho^2/|I|)).
\end{equation*}
Define $R$ by $\pi R^2 = |2B|=4|I|$. Direct computation shows that
\begin{equation*}
\int_{2B}\ln(2\rho)dx=\ln(2\rho)|2B|=4\ln(2\rho)|I|=
2|I|\ln(\rho^2)+2|I|\ln 4
\end{equation*}
and 
\begin{align*}\frac{1}{2\pi}
\int_{2B}\!\ln|x|dx &=\int_0^R r\ln r dr=\frac{1}{2}R^2\ln
R-\frac{1}{4}R^2 \\ &= \frac{|I|}{\pi} \left(\ln
\left(\frac{4|I|}{\pi}\right) - 1\right).\end{align*} So,
$$2\pi||G_m||_{L^1(2B)}=2|I|\ln(\rho^2/|I|)+2|I|(1+\ln(\pi))\leq
2|I|(1+\ln(\pi\rho^2/|I|)).$$

$\hfill\Box$

\section
{\bf Results in dimension 4 or higher}

In this section we discuss uniqueness results in dimensions $n\ge
4$. As in the previous section we give new estimates for unbounded
$D$ to complement the estimates in Lemma~\ref{Global2} for bounded
$D$.

\begin{Thm}\label{n4thm0} Suppose $n\geq 4$. Then \BEL{e4thm1} \max\{c(|D|^{n \over
2(n-1)}+|D|^{2\over n-1}), 1\}||V||_\infty \geq 1 \EE where $c=C (
|k-\pi^2|^{\frac{n-4}{4}}+1)$ and $C$ is a generic constant.
\end{Thm}

\noindent{\it Remark:} From Theorem~\ref{n4thm0}, it is easy to  see
that when $|D|$ is sufficiently small we have the estimate
$||V||_\infty \geq 1$. Note that the estimate \eqref{e4thm1} holds in both cases,
$k\in\mathcal{K}$ and $k\not\in\mathcal{K}$. It is based on the
estimates for $u_m$ in the following two lemmas.

\begin{Lemma}\label{Lprel4}
If $k_m\neq 0$, then
\begin{eqnarray}\label{est_um_n44}
||u_m||_{L^2(I )} &\leq& C(|I|^{n \over
2(n-1)}|k_m|^{\frac{n-4}{2}}+|I|^{2\over n-1}) ||f_m||_{L^2(I )}.
\end{eqnarray}
\end{Lemma}

\begin{Lemma}\label{n4Spec} Suppose that $k=m^2\pi^2$.
 Then
\begin{equation}\label{en4Spec}
\|u_{m}\|_{L^2(I)}\leq C|I|^\frac{2}{n-1}\|f_{m}\|_{L^2(I)}.
\end{equation}
\end{Lemma}

Example \ref{4Ex} and Example \ref{n4Sp*} show that the inequalities
\eqref{est_um_n44} and inequality \eqref{en4Spec} are sharp,
respectively. We prove the above lemmas and the theorem in the next
subsections.

\subsection{Proof of the lemmas}

\noindent {\it Proof of Lemma~\ref{Lprel4}.} We will prove the
following pointwise estimate for the fundamental solution $g_m$
below for $n\geq 4$,
\begin{eqnarray}\label{est-gm4}
|g_m(r)|\leq Cr^{-(n-3)} \left( (|k_m|r)^ {\frac {n-4}2} + 1\right).
\end{eqnarray}
By dilation, it is enough to prove \eqref{est_um_n44} when $|I|=1$.
Let
$$G_m(r)=C r^{-n+3}
 \left((|k_m|r)^{ \frac {n-4}2} +   1 \right)
 $$
 where $C$ is the same as  in \eqref{est-gm4}. Then
 $G_m$ satisfies the requirements of   Lemma~\ref{rearr}, so
  we need to estimate  its $L^1$ norm on the ball  $B_\rho(0) = 2B$ (this defines $\rho$),
  where $B$ is the ball of measure $1$ centered at the origin.
 Thus,
\begin{eqnarray*}
&&||G_m||_{L^1(B(0, \rho))} = |S^{n-2}| \int_0^{ \rho} r^{n-2}
|G_m(r)| dr\\&=& C  \int_0^{ \rho} r (( |k_m|  r)^{ \frac {n-4}2} +
1)dr
  \leq    C  \rho^2
 \left(( |k_m|\rho  )^{\frac{n-4}{2}} +   1 \right)\\
& \leq&  C  \left(|k_m|^{\frac{n-4}{2}}+1\right),
\end{eqnarray*}
which is \eqref{est_um_n44}.

Next, we prove the estimate \eqref{est-gm4}. From \eqref{Eintgm} we
know
\begin{equation}\label{green_n4}
g_m(r)=c(s)(-ik_mr)^{s-\frac{1}{2}}r^{-2s}e^{ik_m r} I_s(-ik_m r)
\end{equation}
where $s=\frac{n-3}{2}\ge \frac 12$ and $ I_s (z) = \int_0^\infty
e^{-t} t^{s-\frac 12}\left(1 +\frac{t}{2z}\right)^{s-\frac 12} dt. $
 Since\newline
$(a+b)^c \le (2a)^c+(2b)^c$ for $a,b,c>0$, we have
\begin{equation}\label{Ineq1}
|I_s (z)|\le  \int_0^\infty\!\!\! e^{-t} t^{s-\frac 12}\left(2^{s-\frac
12} + \left(\frac t{|z|}\right)^{s-\frac 12}\right)dt \leq C (1+
|z|^{-s+\frac 12}).
\end{equation}
Combining \eqref{green_n4} and \eqref{Ineq1} and in view of
$|e^{ik_mr}|\leq 1$, we obtain
\begin{eqnarray*}
|g_m(r)|&\leq& C(|k_m|r)^{s-\frac 12}   r^{-2s} \left|e^{ik_m
r}\right| \left( 1 +
(|k_m|r)^{-s+\frac 12} \right)\\
&\leq &C r^{-2s} \left((|k_m|r)^{s-\frac 12} +1\right) =Cr^{-(n-3)}
\left( (|k_m|r)^ {\frac {n-4}2} + 1\right).\ \Box
 \end{eqnarray*}

\medskip
\noindent {\it Proof of Lemma~\ref{n4Spec}.}  The fundamental
solution of $-\Delta$ in $\mathbb{R}^{n-1}$ for $n\geq 4$ is
\begin{equation*}
g_m=\frac{\Gamma(\frac{n-3}{2})}{4\pi^{\frac{n-1}{2}}}\
\frac{1}{|x|^{n-3}}.
\end{equation*}
Let $G_m=g_m $.  Then $u_m=G_m*f_m$ and $G_m$ satisfies the
requirement in Lemma~\ref{rearr}. We estimate $G_m$ on the ball $2B$
where $B$ is the ball centered at $0$ with measure $|I|$ in
$\R^{n-1}$. Clearly, $2B=B(0,\rho)$ with
$\rho=2(\frac{(n-1)|I|}{|S^{n-2}|})^{\frac{1}{n-1}}$. Direct
computation shows
\begin{align*}
||G_m||_{L^1(B(0, \rho))} &= |S^{n-2}| \int_0^{ \rho} r^{n-2}
|G_m(r)| dr \\ &=C \int_0^{ \rho}r dr=C\rho^2=C|I|^\frac{2}{n-1}.
\end{align*}
By Lemma~\ref{rearr}, Lemma~\ref{n4Spec} follows. $\hfill\Box$

\subsection{Proof of Theorem~\ref{n4thm0}} We prove Theorem~\ref{n4thm0} for both cases $k\not\in\mathcal{K}$
and $k\in\mathcal{K}$.

\medskip

\noindent {\it Proof of Theorem~\ref{n4thm0} when $k\not
\in\mathcal{K}$.} We show that
\begin{equation}\label{enew-n>4}
||u_m||_{L^2(I)}\leq \max\{c(|I|^{n \over 2(n-1)}+|I|^{2\over n-1}),
1\}||f_m||_{L^2(I)}
\end{equation}
for each term $u_m$. Then the inequality \eqref{e4thm1} follows from
\eqref{enew-n>4} and Lemma~\ref{Lmain}. Denote
$m_0=\min\{m\in\mathbb{N}: m^2\pi^2>k\}$. We divide the terms
$u_m$ in three categories.

\noindent (1) For the terms $u_m$ with $m> m_0$, we use the estimate
in   Lemma~\ref{Global1}. In view of $|k_m|^2\geq
3\pi^2>1$, we have
\begin{eqnarray*}
||u_m||_{L^2(I )} \leq \frac{1}{|k_m|^2}||f_m||_{L^2(I )}\leq
||f_m||_{L^2(I )}.
\end{eqnarray*}
So \eqref{enew-n>4} holds for these terms.

\noindent (2) For the term $u_m$ with $m=m_0$, we use different
estimates based on the size of $|k_m|=\delta_+$. If $\delta_+<1$,
the estimate in Lemma~\ref{Lprel4} implies
\begin{eqnarray*}
||u_m||_{L^2(I )} &\leq& C(|I|^{n \over
2(n-1)}|k_m|^{\frac{n-4}{2}}+|I|^{2\over n-1}) ||f_m||_{L^2(I )}\\
&\leq& C(|I|^{n \over 2(n-1)}+|I|^{2\over n-1}) ||f_m||_{L^2(I )}.
\end{eqnarray*}
If $\delta_+\geq 1$, the estimate in Lemma~\ref{Global1} implies
\begin{eqnarray*}
||u_m||_{L^2(I )} \leq \frac{1}{\delta_+^2}||f_m||_{L^2(I )}\leq
||f_m||_{L^2(I )}.
\end{eqnarray*}
So \eqref{enew-n>4} holds for this term.

\noindent (3) For the terms $u_m$ with $m<m_0$. From the estimate in
Lemma~\ref{Lprel4}, in view of $|k_m|\leq \sqrt{k-\pi^2}$, we have
\begin{eqnarray*}
||u_m||_{L^2(I )} \leq C(|I|^{n \over
2(n-1)}|k-\pi^2|^{\frac{n-4}{4}}+|I|^{2\over n-1}) ||f_m||_{L^2(I
)}.
\end{eqnarray*}
So \eqref{enew-n>4} also holds for these terms. Please also note
that if $m_0=1$ (or equivalently $k<\pi^2$), then no term $u_m$
falls in this category, and we just skip this step in this case.
$\hfill\Box$

\medskip

\noindent {\it Proof of Theorem~\ref{n4thm0} when
$k\in\mathcal{K}$.} Since $k\in \mathcal{K}$, we denote
$m_0^2\pi^2=k$. For the term $u_m$ with $m=m_0$, Lemma~\ref{n4Spec}
gives
\begin{equation*}
\|u_m\|_{L^2(I)}\leq C|I|^\frac{2}{n-1}\|f_m\|_{L^2(I)}.
\end{equation*}
For the terms $u_m$ with $m\neq m_0$, similar to the proof of
\eqref{enew-n>4} (see the categories (1) and (3)), we have
\begin{equation*}
||u_m||_{L^2(I)}\leq \max\{c(|I|^{n \over 2(n-1)}+|I|^{2\over n-1}),
1\}||f_m||_{L^2(I)}.
\end{equation*}
In view of $c=C(|k-\pi^2|^{\frac{n-4}{4}}+1)\geq C$,
Lemma~\ref{Lmain} concludes the proof.

$\hfill\Box$

\section{\bf Remarks on sharpness}

Most of our lemmas of the form $||u_m||_2 \le C(\delta,|D|) ||f_m||_2$ are sharp in some sense (the main exceptions are
Lemma \ref{Global2}  and Lemma \ref{Lorenz3}).
We say that an inequality of the form $f(x)\le C g(x)$ is {\it sharp} as $x\to a$ if there is some sequence $x_j\to a$
such that $g(x_j)/f(x_j)$ is bounded.
Usually, $\delta$ or $|D|$ will play the role of $x$ and $a$ will be 0 or $+\infty$.

In general, if a theorem about $V=V_u$ is based on sharp lemmas about $||u_m||_2$, then it is also sharp, so we will focus mainly on our lemmas in this section. We illustrate this principle with a discussion of Theorem \ref{ThmX}, which concludes with
$$ \max\left\{\frac 1 {\delta_+^{2}},\ \frac{ C\rho}{\delta_-}\right\} ||V||_\infty >1. $$
The formula $\delta_+^{-2}$ in it comes from  the inequality $||u_m||_2  \leq \delta_+^{-2}||f_m||_2  \leq \delta_+^{-2}||V_{u_m}||_\infty ||u_m||_2   $ in Lemma \ref{Global1}. Example \ref{forGlobal} shows that this inequality  is sharp.
To examine the sharpness of $\delta_+^{-2}$ in Theorem \ref{ThmX}, we set $u(x,y)=  u_m(x)\sin(\pi m y)$ (so, the Fourier series of $u$ has only one term) and will show that $u$ satisfies \eqref{eqk} with fairly small $||V||_{L^\infty(D)}$. We note that $u_m = g_m$ off $I$ (by construction in Example \ref{forGlobal} and other similar examples) and it
satisfies the equation \eqref{umPDE2} with $f_m= u_m V_m$ on $S$, with $V_ m = -\Delta u_m/u_m -k_m^2$ supported on $I$.
So, for $(x,y)\not\in D$,
\begin{eqnarray*} (-\Delta -k)u &=& (-\Delta -k)g_{m}(x)\sin(m\pi y) \\ &=& \left( -\Delta_x  - k_{m} ^2 \right)g_m(x) \sin(m\pi y) =0. \end{eqnarray*}
On $D$, observe that $$ -\Delta u = [\frac{-\Delta_x v}{v} + (\pi m)^2] u = [V_m+  k_{m} ^2 + (\pi m)^2] u= [V_m+ k] u.$$
So, $u$ satisfies \eqref{eqk} with  $V(x,y) = V_m(x)$ supported on $D$, and \newline $||V||_{L^\infty(D)}$ $ = ||V_m||_{L^\infty(I)} \approx \delta_+^{2}$
as $\delta_+ \to \infty$. Thus, the formula $\delta_+^{-2}$ in Theorem \ref{ThmX} is sharp in this sense. Note that $ \max\left\{\frac 1 {\delta_+^{2}},\ \frac{ C\rho}{\delta_-}\right\} =
\frac 1 {\delta_+^{2}}$ can occur also if $\delta_+$ is large.

However, we do not claim that Lemma \ref{Global2}  is sharp in all dimensions. So, we cannot claim that the term
$\frac{ C\rho}{\delta_-}$ in Theorem \ref{ThmX} is sharp, for example.
We will leave this kind of reasoning to the reader for the other theorems in the paper and will focus here on the sharpness of the lemmas.

One advantage of this approach is that the lemmas are dilation invariant (most of the theorems on $u$ are not).
The reader may note that since Example \ref{forGlobal} is based on Lemma \ref{dLarge}, its potential $V_m$ has support on $I=B_1(0)$.
But this is not necessary. Since Lemma \ref{Global1}  is dilation invariant, Example \ref{forGlobal} can be easily revised so that $I=B_R(0)$ for any given $R>0$. Similar remarks apply to all examples in this section which concern $||V_m||_\infty$ or $||u_m||_2$, but not to
the few examples about $||V||_\infty$, such as Example \ref{n2SpV}. So, we can focus on $u_m$ and sharpness in terms of $k$ (or $\delta$).

The first subsection below contains some general purpose lemmas and
several examples which are direct consequences of those. The next
one contains special constructions needed for $n=2$ and $n=3$. The
last one deals with the cases $k< \pi^2$ and $k\in \mathcal{K}$.

\subsection{Patching lemmas}
We use Lemma \ref{parab} below to construct examples with fairly small potentials
based on radial parabolic interpolation near the origin. We are
given $k$ and $m$ and some radial function $v(r)$ (such as $g_m(r)$) defined for $r=|x|>{1}$
in $\R^{n-1}$ that satisfies $-\Delta_x v = k_m^2 v$.
Based on the boundary values of $v$ where $r={1}$, we want to extend $v$ to $I = B_{{1}}(0) \subset \R^{n-1}$ such that
$-\Delta_x v = V_v v +k_m^2 v$, with $||V_v||_{L^\infty(B_{1}(0))}$ fairly small.
Since the extension involves simple quadratic interpolation, we refer to this as {\it parabolic patching}.
In the Lemma below and the examples that follow, we  will let $V_v= -\Delta_x v/v-k_m^2$. When $v=u_m$, we will write $V_m$ instead of $V_{u_m}$.

The first lemma is intended for the real-valued $v$ that typically occur when $k_m^2 = k-\pi^2m^2 <0$.

\begin{Lemma}\label{parab}
Suppose $v$ is given, for $r\ge {1}$, with a well-defined $v'({1})\le 0$.
Let $v(r) = a-br^2 > 0$ for $r\le {1}$, where $-v'({1}) = 2b$. Then
$|V_v +k_m^2|  \le \frac{n-1}{{1}}|v'({1} )/v({1})|$ on $I = B_{{1}}(0)\subset \R^{n-1}$.
\end{Lemma}
\noindent{\it Proof:} $v'(r) = -2br$ and $v''(r) = -2b$ and (on
$\R^{n-1}$), we have $\Delta_x v (r)= -2b(n-1) =
\frac{n-1}{{1}}v'({1})$. $\hfill\Box$
\medskip
Lemma \ref{parabC}  below is similar, but is intended for $\pi^2m^2 <k$, when the functions involved will be complex-valued and
$k_m\in \R$. It is most useful when $k_m$ is small.

\begin{Lemma}\label{parabC}  Let $\psi(r) = A- B r^2$, where $A$, $B\in\C$ and $0\leq r\leq {1}$.  Assume that  $\Re (B) \ge 0$ and   $\Re(\psi({1}))>0$.
Then, on $I = B_{{1}}(0)\subset \R^{n-1}$,
 $$|V_\psi +k_m^2|   \le \frac{(n-1)|\psi'({1})|}{{1} \Re(\psi ({1}))}.$$
\end{Lemma}
\noindent{\it Proof:} It is similar to Lemma \ref{parab}. $\Delta
\psi(r)$ simplifies to $(n-1)(- 2B)= \frac{(n-1) \psi'({1}) }{{1}  }
$ and $|\psi (r)| \ge \Re (\psi (r)) \ge \Re (\psi ({1})) $.
$\hfill\Box$
\medskip

Lemma \ref{dLarge} combines patching with Lemma \ref{Lg1/g}.
We will use it to prove that several lemmas from Sections 2 to 5 are sharp.
\begin{Lemma}\label{dLarge}
If $m^2\pi^2 \ne k$, there exists a nontrivial solution  $u_m$ of \eqref{umPDE} with $V_{m}$ supported in $B_1(0)$
such that
\begin{itemize}\item[1)]
$||V_{m} + k_m ^2||_\infty\le c_n   (|k_m|+1)$ if $n\ge 3$, and

\item[2)]
$||V_{m} + k_m ^2||_\infty\le c_n|k_m|$ if $n=2$.
\end{itemize}
\end{Lemma}

\noindent {\it  Proof:} Let $u_m(r) = g_m(r) $ for $|r|>1$. Extend
this to $|r|\leq 1$ using Lemma \ref{parab} when $m > k/\sqrt \pi$,
or  Lemma \ref{parabC} when $m < k/\sqrt \pi$. By Lemma \ref{Lg1/g},
$$||V_{m} +k_m ^2||_{L^\infty (D)} \le (n-1) \left|\frac{g'_m(1)}{g_m(1)}\right|  \leq
 c_n (1+ |k_m|)$$ when $n\ge 3$. The case $n=2$ also follows from Lemma \ref{Lg1/g}. $\hfill\Box$
\medskip

\begin{Ex}\label{forGlobal}  (Lemma \ref{Global1} is sharp.) Suppose  $n\ge 2$ and $m>\sqrt k/\pi$.  Then
  $\lim_{\delta_+\to\infty} ||V_{m}||_{\infty} /\delta_+^{2}=1$ can occur.
\end{Ex}

Given $\delta_+$, we can create an admissible solution using Lemma
\ref{dLarge} with $|k_m|=\delta_+$ and $||V_m+k_m^2|| <
c_n(1+|k_m|)$. With this and Lemma \ref{Global2}, we have $|k_m|^2
\leq ||V_{m}||_{\infty} \leq  |k_m|^2+ c_n(|k_m|+1)$. $\hfill\Box$

\begin{Ex}\label{2Ex}  (Lemma \ref{n2lemma} is sharp.) Suppose $n=2$, and $m\ne \sqrt k/\pi$.
For all $\delta \le 1$, $||V_{m}||_{\infty}  \le C\delta$ can occur.
\end{Ex}

That is, there is an absolute constant $C$  such that, for a  given $\delta \le 1$, we can find an admissible solution $u$ such that
$||V_{m}||_{\infty}  \le C\delta$.

\begin{Ex}\label{4Ex}
(Lemma \ref{Lprel4} is sharp.) Suppose $n\ge 4$ and $m\ne \sqrt
k/\pi$. For all $\delta \le 1$, $||V_{m}||_{\infty}  \le C $ can
occur.
\end{Ex}

The constructions for Examples \ref{2Ex} and \ref{4Ex} are identical
to Example \ref{forGlobal} above. $\hfill\Box$

\subsection{Special examples for n=2 and 3}

In this subsection we discuss two cases which are not covered by
Lemma 6.3.

\begin{Ex}\label{n2mD} (Lemma~\ref{Global2} is sharp when $n=2$; see also Lemma
\ref{n2lemma}.) Suppose  $n=2$, $I=[-1,1]$ and $m<\sqrt k/\pi$. For all $\delta_- >1$, $||V_m||_\infty \le C \delta_-$ can occur.
\end{Ex}

Note that Lemma \ref{dLarge} is useless here, because  when
$\delta_-$ is large, we cannot expect $||V_{m}+\delta_-^2||_\infty
\approx ||V_{m}||_\infty$. Instead, choose $B \approx 1$  so  that
$|\cos(\delta_- B)| = \frac 12$. We construct an admissible $u_m$
that solves $-u_m''(x) = \delta_-^2 u_m$ off the set $I = [-B,B]
\approx [-1,1]$, with $||V_m||_\infty \approx \delta_-$.  Set
\BEL{hardEx}u_m(x) = \phi(x)\sin(|\delta_- x|) -i\cos(\delta_-
x).\EE We  construct an even function $\phi(x)$ on $\R^1$ such that
$\phi(x)\equiv 1$ for $|x|\ge B$. So, $iu_m = e^{i\delta_- |x|}$
there, which is admissible and satisfies $-u_m''(x) = \delta_-^2
u_m$. We must define  $\phi(x)$ on $[0, B]$. On each maximal
interval $I_j\subset [0, B]$ where $|\cos(\delta_- x)| \le \frac
12$, $\phi$ will be a constant $c_j \ge 0$ to be defined
recursively.

Let $J_j=(a ,b )$ be a maximal interval where $|\cos(\delta_- x)| >
\frac 12$, and assume $\phi(a ) = c_j$ has already been defined
(with $c_0 = \phi(0)=0$). Note that $b -a  \approx  1/\delta_ -$.
Let $d $ be the midpoint of $J_j$ and define $\phi(x) = c_j +
\delta_-(x-a )^2$ on $(a ,d ]$. Define $\phi(x-d ) -\phi(d) =
\phi(d) - \phi(2d-x)$ on $(d,b]$ (giving local anti-symmetry around
the point $(d,\phi(d))$) and let $c_{j+1} =\phi(b)$. Note that $\phi
\in C^1$ is nondecreasing and $|\phi''(x)| \le \delta_-$. The
maximum value of $\phi'(x)$ occurs at each $d$ and is $2\delta_-
(d-a) \approx 1$. Recall $V_{m} = -u_m''/u_m - \delta_-^2$. From
\rq{hardEx}, $$u_m'' = \phi'' s +2 \phi' s' + \phi s'' + \delta_-^2
i\cos(\delta_- x)$$ where $s(x) = \sin(\delta x)$. Note $s'' =
-\delta_-^2 s$ so that $u_m'' = -\delta_-^2 u_m +E$ with $E = \phi''
s +2 \phi' s'$. By construction, $E=0$ unless $|\cos(\delta_- x)|
\ge \frac 12$. For such $x$, $|E|\le C \delta_- \le C \delta_-
|u_m|$.  This shows $|V_{m}| \le |E/u_m| \le C \delta_-$.
$\hfill\Box$

 \medskip  This example concludes the proof of sharpness for all  Lemmas in dimension $n=2$; we now turn to $n=3$.

\begin{Ex}\label{Forn3S} (Lemma \ref{n3S} is sharp.)
Suppose  $n=3$ and $m\ne\sqrt k/\pi$. Then for all sufficiently small $\delta$,
$||V_m||_\infty |\log(\delta)|\le C$
can occur.
 \end{Ex}
We will assume that  $m <\sqrt k/\pi$; the case $m >  \sqrt k/\pi$
is similar. Let $|k_m|=\delta  <1$. Let $u_m(r) = 4 g_m(r)=
-Y_0(\delta r)+ iJ_0(\delta r)$ for $r\ge 0.5$. Extend $u_m$ to
$R^2$ by parabolic patching as in Lemma \ref{parabC}. We need some
estimates on the Bessel functions $Y_0$ and $J_0$, to apply a fairly
obvious modification of this lemma (with 1 replaced by 0.5) to
$u_m$. For the basic properties of these functions, see the classic
\cite{Er}. Using the identity $J_0'(x)= J_{-1}(x)= -J_{1}(x)$ and
the Poisson representation formula, we see that $J_0(x) \approx 1$
and $|J_0'(x)| \le 1$  in $[0, 0.5]$. Also recall that  $Y_0(x)= c
J_0 (x) |\ln(x/2)+\gamma|+b(x)$ where $b$ is an analytic function,
and $\gamma \approx 0.577$ is  the Euler-Mascheroni constant. Since
$b$ is defined as an alternating power series, we see that $0\le
b(x) \le 1/16$ and $b'(x)\le 2$ for $0\leq x\leq 0.5$. When $\delta$
is sufficiently small,
$$\Re \, u_m (0.5) = -Y_0(\delta /2) =  cJ_0(\delta /2)|\ln(\delta /4)+\gamma| +b(\delta /2) \ge c |\ln(\delta)| ,$$
and $|{\rm Re}\ u_m' (0.5)| = |\delta Y_0'(\delta /2)| \le \delta C
|\ln(\delta)| \le C .$ Also, $|{\rm Im}\  u_m'(0.5)| = |\delta
J_0'(\delta /2)| \le C$. Noting that  $\delta^2 \le
C/|\ln(\delta/2)|$, the modified Lemma \ref{parabC} implies $||V_m||
\le |k_m|^2 + C/|\ln(\delta/2)| \le C/|\ln(\delta/2)|$ as desired.
$\hfill\Box$

\subsection{Sharpness proofs for  $k<\pi^2$ and $k\in \mathcal K$.}
In this subsection, we study the sharpness results for some special
values of $k$, in all dimensions. In contrast to the previous
examples in this section, we will focus on $||u_m||_2$ or $||V||$,
rather than $||V_m||$.
\medskip
\begin{Ex}\label{k<0} (Corollary \ref{Csmallk} is sharp.)
Assume $n\ge 2$ and $k<\pi^2$. Then, for every $\epsilon>0$,  $||V||_\infty \le
\delta_+^2+\epsilon$ can occur.
\end{Ex}

Note that $\delta_+ = k_1 =\sqrt{\pi^2-k}$. Let $r_0$ be a large
positive number to be specified later. Let $s(r) = e^{-\delta_+ r}$
for $r> r_0$ and let $ s(r) =A r^2 +B$ for $r\le r_0$, where $A =
-\delta_+e^{-\delta_+ r_0}/(2r_0)$ and $B = (1+
\delta_+r_0/2)e^{-\delta_+ r_0}$ are chosen so to make $s$
differentiable. Let $u(x,y)=s(|x|)\sin(\pi y)$ on $S$. A direct
computation shows that $-\Delta u= (V+\pi^2-\delta_+^2)u = (V+k)u$,
where
$$V(r,y)= \begin{cases} \frac{\delta_+(n-2)}{r} \ \mbox{ if} \ r>r_0
\\
{ \delta_+^2+ \frac{2 \ \delta_+ (n-1)}{r_0^2 \delta_++2
r_0-\delta_+r^2}}\ \mbox{ if} \ r\leq r_0. \end{cases}
$$
So, $ ||V||_\infty = \delta_+^2+  \frac{\delta_+(n-1)}{r_0} <
\delta_+^2 +\epsilon$ for  $r_0$ sufficiently large. $\hfill\Box$

\medskip

Our next two examples are for $n=2$ with $k\in \mathcal K$. It seems unlikely that examples for $||V_{m}||_\infty$ exist,
but we can still show our results are sharp.

\begin{Ex}\label{n2Sp*} (Lemma \ref{n2Spec} is sharp.)
Given $\rho>0$ and $|I|\le 2\rho$, $\|u_m\|_{L^2(I)}\ge C\rho |I|\
\|f_m\|_{L^2(I)}$ can occur.
\end{Ex}
By dilation (see the remarks below \rq{diln})  we can assume that
$\rho =1$. Suppose $|I|\le 2$ is given. Let $u_m = 2/3\ - |x|$, with
these modifications: a) let $u_m \equiv 0$ off of $(-1, 1)$; b)
smooth $u_m$ to be in $C^1(\R)$, so that $f_m = -u_m'' \equiv 0$
except on intervals of length $|I|/3$ centered at $-2/3,\ 0$ and
$2/3$, where $f_m$ is constant on each of these intervals. Let $I
\subseteq [-1,1]$ be the union of these 3 intervals. Since $f_m
\approx |I|^{-1}$ on $I$, $||f_m||_2 \approx |I|^{-\frac 12}$ and
$||u_m||_2 \approx |I|^{\frac 12}$, as desired.$\hfill\Box$
\medskip
\begin{Ex}\label{n2SpV} (Theorem \ref {n2main} part (2) is sharp when $\rho=1$.)
Let $n=2$ and $|D|\le 2\rho$. Then $||V||_\infty   \le   C|D|^{-1}$
can occur.

\end{Ex}
For simplicity, let $k = 4\pi^2$. We may assume $|D|<0.1$. Let $u_1
= 2\cos(\sqrt{3}\pi x)$, and define $u_2 \approx 0.1 -|x|$, smoothed
out as in Example \ref{n2Sp*} so that $u_2''$ is supported on a set
of measure $|D|$. Note that $0\le 2u_2 \le u_1$ on the support of
$u_2$. Let $u(x,y)= u_1(x) \sin  (\pi y)+ u_2(x) \sin(2\pi y)$,
which satisfies \eqref{eqk} with $V= - \frac{2\cos( \pi y)
u_2^{\prime\prime }(x)} {u_1(x)  + 2 u_2(x) \cos( \pi y)} $. Let $I$
be  the support of $u_2''$. If $V(x)\not = 0$, then $x\in I$ and
$u_1(x)  + 2 u_2(x) \cos( \pi y) \approx u_1(x) \approx 1$. So, $|V|
\le C|u_2''(x)| \approx |D|^{-1}$, as in Ex \ref{n2Sp*}.
$\hfill\Box$
\medskip
\begin{Ex}\label{n3Sp*} (Lemma \ref{n3Spec} is sharp)  Let $n=3$  and $k=m\pi^2$. If
$\rho>0$ and $|I|\le \pi\rho^2$,  then $||u_m||_2 \ge C |I|\ln(1+ \pi \rho^2 / |I|) ||f_m||_2$ can occur.
\end{Ex}
Set $u = u_m(x) \sin(m\pi y)$ with $k=m^2\pi^2$, so that $-\Delta_x u_m = f_m$.
By dilation (see \rq{diln}), we may assume $\rho=1$.
We may assume $|I|$ is small,
so that $\ln(1+ \pi \rho^2 / |I|)\approx -\ln(|I|)$.
Let $v_1(r) = -\ln(r)$ for $r<1$ and zero otherwise.
Define $A_{\alpha, \beta} = \{x\in \R^2: \alpha \le |x| \le \beta\}$.
Let $I_1 = A_{0,a}$, such that $|I_1| =|I|/2$ and let
$I_2 = A_{b,1}$, such that $|I_2| =|I_1|$. Let $v_2=v_1$ off $I=I_1\cup I_2$.
Define $v_2 \in C^1(I_1)$ by: $v_2(a)=v_1(a)$ and $v_2'(a)=v_1'(a)$;
and on $(0,a)$ let $-\Delta v_2  = c_1$  be a constant such that $v_2'(0)=0$.
By Green's identity,
$$c_1 |I_1| = - \int_{I_1} \Delta v_2 \ dx = - 2\pi a v_2'(a) = 2\pi.$$
So, $||\Delta v_2||_{L^2(I_1)} = c_1 |I_1|^{\frac 12} = C |I_1|^{-\frac 12}$.
Define $v_2 \in C^1(I_2)$ by:
$v_2(b)=v_1(b)$ and $v_2'(b)=v_1'(b)$; and on $(b,1)$
let $-\Delta v_2 = c_2$ be a constant such that $v_2'(1)=0$.
Let $\epsilon = v_2(1) >0$ and define $u_m = v_2(r) -\epsilon$. So,
$u_m \in C^1$ and $u_m(1)=0$.
Since $|I|$ can be assumed to be small, we may also assume $\epsilon < 1/10$. As above
$$c_2 |I_2| = - \int_{I_2} \Delta u_m \ dx = - 2\pi b v_2'(b) = 2\pi.$$ So,
$||\Delta u_m||_{L^2(I_2)} = C |I_2|^{-\frac 12}$, as for $I_1$. On
$I_1$, $u_m(r) \ge v_1(b) - \epsilon \approx -\ln(|I|)$, since $a^2
\approx |I|$ and $\epsilon$ is negligible. So $||u_m||_{L^2(I)} >
||u_m||_{L^2(I_1)} > -C\ln(|I|) |I_1|^{\frac 12} \ge -C\ln(|I|) |I|
||\Delta u_m||_{L^2(I)} $. $\hfill\Box$

\medskip
\begin{Ex}\label{n4Sp*} (Lemma \ref{n4Spec} is sharp). Assume $n\ge 4$ and $|I|$ is given.
Then $||u_m||_2  \ge  C  |I|^{\frac 2{n-1}} ||f_m||_2$ can occur.
\end{Ex}
Define $u$, $k$ and $m$ as in Ex. \ref{n3Sp*}. Let $u_m(r)= r^{3-n}$
for $r > r_0$, where $|B_{r_0}|=|I|$ (so $c_n r_0^{n-1} = |I|$). On
$I = B_{r_0}$ define $u_m$ so that $u_m\in C^1(\R^{n-1})$ and $f_m =
-\Delta u_m$ is a constant. Again $||f_m||_1 \approx 1$ and $f_m
\approx |I|^{-1}$ and $||f_m||_2 \approx |I|^{-\frac 12}$. Again
$||u_m||_2 > |I|^{\frac 12}$ min $u_m|_I = |I|^{\frac 12} r_0^{3-n}
\approx |I|^{\frac 12} |I|^{\frac {3-n}{n-1}}$. So,
$||u_m||_2/||f_m||_2 \ge C_n |I|^{\frac 2{n-1}}$. $\hfill\Box$

\end{document}